\documentclass{article}
\usepackage{amsmath}
\usepackage{amssymb}
\usepackage{amsfonts}

\newtheorem{theorem}{Theorem}[section]

\newtheorem{definition}[theorem]{Definition}

\newtheorem{lemma}[theorem]{Lemma}

\newtheorem{remark}[theorem]{Remark}

\newenvironment{proof}[1][Proof]{\textbf{#1.}\;}{\vspace{-0,25 in} \begin{flushright}\rule{0.3em}{0.3em}
\end{flushright}}
\def\qed{\ifhmode\textqed\fi
   \ifmmode\ifinner\quad\qedsymbol\else\dispqed\fi\fi}
\def\textqed{\unskip\nobreak\penalty50
    \hskip2em\hbox{}\nobreak\hfil\qedsymbol
    \parfillskip=0pt \finalhyphendemerits=0}
\def\dispqed{\rlap{\qquad\qedsymbol}}

\begin{document}

\title{Quasi-Nilpotent Operators on Locally Convex Spaces \thanks{%
This paper was written during my visit to Departamento de Analisis
Matematico, from Facultad de Ciencias Matematicas de la Universidad
Complutense de Madrid, and was supported by the MEdC-ANCS CEEX grant
ET65/2005 contract no. 2987/11.10.2005 and the M.Ed.C. grant C.N.B.S.S. contract no. 5800/09.10.2006. }}
\author{Sorin Mirel Stoian}
\date{}
\maketitle

\begin{abstract}
In this article we extend the notion of quasi-nilpotent equivalent operators, introduced by Colojoara and Foias \cite{co1} for Banach spaces, to the class of bounded operators on sequentially complete locally convex spaces. 
\end{abstract}

{\small AMS 2000 Subject Classification: 47B40, 47B99}

\section{ Introduction}

The class of quasi-nilpotent equivalent operators on a Banach space was
introduced by Colojoara and Foias \cite{co1}. The aim of this paper is to
search if we can extend this theory to the class of bounded operators on
sequentially complete locally convex spaces.

Any family $\mathcal{P}$ of seminorms which generate the topology of a locally
convex space $X$ (in the sense that the topology of $X$ is the coarsest with
respect to which all seminorms of $\mathcal{P}$ are continuous) will be
called a calibration on $X$. The set of all calibrations for $X$ is denoted
by $\mathcal{C}(X)$ and the set of all principal calibration by $\mathcal{C}%
_{0}(X)$.

An operator $T$ on a locally convex space $X$ is quotient bounded with
respect to a calibration $\mathcal{P}\in \mathcal{C}(X)$ if for every
seminorm $p\in \mathcal{P}$ there exists some $c_{p}>0$ such that 
\begin{equation*}
p\left( Tx\right) \leq c_{p}p\left( x\right) ,\left( \forall \right) x\in X.
\end{equation*}

The class of quotient bounded operators with respect to a calibration $%
\mathcal{P}\in \mathcal{C}(X)$ is denoted by $Q_{\mathcal{P}}(X)$. For every 
$p\in \mathcal{P}$ the application $\hat{p}:Q_{\mathcal{P}}(X)\rightarrow \mathbf{R}$
defined by 
\begin{equation*}
\hat{p}(T)=\inf \{\;r>0\;|\;p(Tx)\leq rp\left( x\right) ,\left( \forall
\right) x\in X\},
\end{equation*}%
is a submultiplicative seminorm on $Q_{\mathcal{P}}(X)$, satisfying the
relation $\hat{p}(I)=1$, and has the following properties

\begin{enumerate}
\item $\hat{p}(T)$=$\sup\limits_{p\left( x\right) =1}p\left( Tx\right)
=\sup\limits_{p\left( x\right) \leq 1}p\left( Tx\right) $, $\left( \forall
\right) p\in \mathcal{P}$, $\left( \forall \right) q\in \mathcal{Q}$;

\item $p\left( Tx\right) \leq \hat{p}\left( T\right) p\left( x\right) $, $%
\left( \forall \right) x\in X$.
\end{enumerate}

We denote by $\hat{{\mathcal{P}}}\ $ the family $\{\hat{p}\;|\;p\in \mathcal{%
P}\}$. If $T\in Q_{\mathcal{P}}(X)$ we said that $\alpha \in \mathbb{C}$ is
in the resolvent set $\rho (Q_{\mathcal{P}},T)$ if there exists $(\alpha
I-T)^{-1}\in Q_{\mathcal{P}}(X)$. The spectral set $\sigma (Q_{\mathcal{P}%
},T)$ will be the complement set of $\rho (Q_{\mathcal{P}},T)$.

An operator $T\in Q_{\mathcal{P}}(X)$ is a bounded element of the algebra $%
Q_{\mathcal{P}}(X)$ if it is bounded element in the sense of G.R.Allan \cite%
{al}, i.e some scalar multiple of it generates a bounded semigroup. The
class of the bounded elements of $Q_{\mathcal{P}}(X)$ is denoted by $(Q_{%
\mathcal{P}}(X))_{0}$. If $r_{\mathcal{P}}(T)$ is the radius of boundness of the
operator $T$ in $Q_{\mathcal{P}}(X)$, i.e. 
\begin{equation*}
r_{\mathcal{P}}(T)=\inf \{\alpha >0\;\,\mid \;\,\alpha ^{-1}T\text{
generates a bounded semigroup in }Q_{\mathcal{P}}(X)\},
\end{equation*}%
then in \cite{al} is proved that 
\begin{equation*}
r_{\mathcal{P}}(T)=\sup \{\;\,\underset{n\rightarrow \infty }{\limsup }%
\left( \hat{p}\left( T^{n}\right) \right) ^{1/n}\mid \;p\in \mathcal{P}\}.
\end{equation*}

The Waelbroeck resolvent set $\rho _{W}(Q_{\mathcal{P}},T)$ of an operator $%
T\in (Q_{\mathcal{P}}(X))_{0}$ is the subset of elements of $\lambda _{0}\in 
\mathbb{C}_{\infty }=\mathbb{C}\cup \left\{ \infty \right\} $, for which
there exists a neighborhood $V\in \mathcal{V}_{(\lambda _{0})}$ such that:

\begin{enumerate}
\item the operator $\lambda I-T$ is invertible in $Q_{\mathcal{P}}(X)$ for
all $\lambda \in V\backslash \{\infty \}$

\item the set $\{~\left( ~\lambda I-T~\right) ^{-1}|~\lambda \in V\backslash
\{\infty \}~\}$ is bounded in $Q_{\mathcal{P}}(X)$.
\end{enumerate}

The Waelbroeck spectrum of $T$, denoted by $\sigma _{W}(Q_{\mathcal{P}},T)$,
is the complement of the set $\rho _{W}(Q_{\mathcal{P}},T)$ in $\mathbb{C}%
_{\infty }$. It is obvious that $\sigma (Q_{\mathcal{P}},T)\subset \sigma
_{W}(Q_{\mathcal{P}},T)$. An operator $T\in Q_{\mathcal{P}}(X)$ is regular
if $\ \infty \notin \sigma _{W}(Q_{\mathcal{P}},T)$, i.e. there exists some $%
t>0$ such that:

\begin{enumerate}
\item the operator $\lambda I-T$ is invertible in $Q_{\mathcal{P}}(X)$, for
all $\mid \lambda \mid >t$

\item the set $\left\{ R\left( \lambda ,T\right) \mid \mid \lambda \mid
>t\right\} $ is bounded in $Q_{\mathcal{P}}(X)$.
\end{enumerate}

Given $(X,\mathcal{P}$) a locally convex space, for each $p\in \mathcal{P}$ we
denote by $N^{p}$ the null space and by $X^{p}$ the quotient space $X/N^{p}$%
. For each $p\in \mathcal{P}$ consider the canonical quotient map $\pi
_{p}:X\rightarrow X/N^{p}$ given by relation 
\begin{equation*}
\pi _{p}(x)=x_{p}\equiv x+N^{p},\left( \forall \right) x\in X,
\end{equation*}%
(from $X$ to $X^{p}$) which is an onto morphism. It is obvious that $X_{p}$
is a normed space, for each $p\in \mathcal{P}$, with norm $||\bullet ||_{p}$
defined by 
\begin{equation*}
\left\Vert x_{p}\right\Vert _{p}=p\left( x\right) ,\left( \forall \right)
x\in X.
\end{equation*}

Consider the algebra homomorphism $T\rightarrow T^{p}$ of $Q_{\mathcal{P}%
}(X) $ into $\mathcal{L}(X^{p})$ defined by 
\begin{equation*}
T^{p}\left( x_{p}\right) =(Tx)_{p},\left( \forall \right) x\in X.
\end{equation*}

This operators are well defined because $T(N^{p})\subset N^{p}$. Moreover,
for each $p\in \mathcal{P}$, $\mathcal{L}(X_{p})$ is a unital normed algebra
and we have 
\begin{equation*}
\left\Vert T^{p}\right\Vert _{p}=\sup \left\{ \left\Vert
T^{p}x_{p}\right\Vert _{p}\;\left\vert \;\left\Vert x_{p}\right\Vert
_{p}\leq 1\text{ for }x_{p}\in X_{p}\right. \right\}
\end{equation*}%
\begin{equation*}
=\sup \left\{ p\left( Tx\right) \;\left\vert \;p\left( x\right) \leq 1\text{
for }x\in X\right. \right\}=\hat{p}(T)
\end{equation*}

For every $p\in \mathcal{P}$ consider the normed space $(\tilde{X}%
^{p},\left\Vert {\bullet }\right\Vert _{p})$ the completition of $%
(X_{p},\left\Vert {\bullet }\right\Vert _{p})$. If $T\in Q_{\mathcal{P}}(X)$%
, then the operator $T^{p}$ has an unique continuous linear extension $\tilde{%
T}^{p}$ on $(\tilde{X}^{p},\left\Vert {\bullet }\right\Vert _{p})$ and 
\begin{equation*}
\sigma (Q_{\mathcal{P}},T)=\underset{p\in \mathcal{P}}{\cup }\sigma (\tilde{T%
}_{p})=\underset{p\in \mathcal{P}}{\cup }\sigma (T_{p}).
\end{equation*}

\section{Bounded Operators with SVEP}

\begin{lemma}
{\label{lemma:sv1}} If $(X,\mathcal{P})$ is a sequentially complete locally
convex space and $T\in (Q_{\mathcal{P}}(X))_{0}$, then%
\begin{equation*}
\overset{\circ }{\rho }(Q_{\mathcal{P}},T)=\rho _{W}(Q_{\mathcal{P}},T).
\end{equation*}
\end{lemma}

\begin{proof}
Assume that there exists $\lambda _{0}\in \rho (Q_{\mathcal{P}},T)\backslash
\rho _{W}(Q_{\mathcal{P}},T)$ such that $\lambda _{0}\in \overset{\circ }{%
\rho }(Q_{\mathcal{P}},T)$. Since $\lambda _{0}\notin \rho _{W}(Q_{\mathcal{P%
}},T)$, then for each neighborhood U of $\lambda _{0}$ the set 
\begin{equation*}
\{~\left( ~\lambda I-T~\right) ^{-1}|~\lambda \in U~\}
\end{equation*}%
is not bounded in $Q_{\mathcal{P}}(X)$. Let $U\in \rho (Q_{\mathcal{P}},T)$
an open set such that $\lambda _{0}\in U$. This implies that there exists $%
\lambda _{1}\in U$ and $p\in \mathcal{P}$ such that for every $n\in N$ there
exists $x_{n}\in X$ ($p(x_{n})\neq 0$) with the property%
\begin{equation*}
p\left( R(\lambda _{1},T)x_{n}\right) >np\left( x_{n}\right) ,
\end{equation*}

Therefore, for $y_{n}=R(\lambda _{1},T)x_{n}$ we have%
\begin{equation*}
p\left( y_{n}\right) >np\left( (\lambda _{1}I-T)y_{n}\right) ,
\end{equation*}%
which implies that $\lambda _{1}\in \sigma _{a}(Q_{\mathcal{P}},T)\subset
\sigma (Q_{\mathcal{P}},T)$ (see \cite{k2}). This contradicts the
supposition we made, so lemma is proved.
\end{proof}

\begin{definition}
If $(X,\mathcal{P})$ is a sequentially complete locally convex space we say
that the operator $T\in (Q_{\mathcal{P}}(X))_{0}$ has the single-valued \
extension property (we will write SVEP) if for any analytic function $%
f:D_{f}\rightarrow X$, where $D_{f\text{ }}\subset \mathbb{C}$ is an open
set, with the property 
\begin{equation*}
(\lambda I-T)f(\lambda )\equiv 0_{X},\left( \forall \right) \lambda \in
D_{f},
\end{equation*}%
results that $f\equiv 0,\left( \forall \right) \lambda \in D_{f}$.
\end{definition}

\begin{definition}
Let $(X,\mathcal{P})$ be a sequentially complete locally convex space and $%
T\in (Q_{\mathcal{P}}(X))_{0}$. For every $x\in X$ we say that the analytic
function $f_{x}:D_{x}\rightarrow X$ is an analytic extension of the function 
$\lambda \rightarrow R(\lambda ,T)$ if $D_{x}$ is an open set such that $%
\rho _{W}(Q_{\mathcal{P}},T)\subset D_{x}$ and 
\begin{equation*}
(\lambda I-T)f(\lambda )\equiv x,\left( \forall \right) \lambda \in D_{x}.
\end{equation*}%
Denote by \ $\rho _{T}(x)$ the set of all complex number $\lambda _{0}$ for
which there exists an open set $D_{\lambda _{0}}$, such that $\lambda
_{0}\in D_{\lambda _{0}}$, and an analytic function $f_{x}:D_{\lambda
_{0}}\rightarrow X$ which has the property 
\begin{equation*}
(\lambda I-T)f_{x}(\lambda )\equiv x,\left( \forall \right) \lambda \in
D_{x}.
\end{equation*}%
The set $\sigma _{T}(x)$ will be the complement of the set $\rho _{T}(x)$.
\end{definition}

\begin{remark}
\begin{enumerate}
\item In the case of bounded operators on a Banach space we have the
condition $\rho (T)\subset D_{x}$, but the lemma \ref{lemma:sv1} implies
that this conditions in the case of quotient bounded operators on
sequentially complete locally convex space is naturally replaced by the
condition $\rho _{W}(Q_{\mathcal{P}},T)\subset D_{x}$.

\item It is known that for a locally bounded operator $T\in Q_{\mathcal{P}%
}(X)$ we have the equalities $\ $ 
\begin{equation*}
\rho (Q_{\mathcal{P}},T)=\rho _{W}(Q_{\mathcal{P}},T)=\rho (T),
\end{equation*}%
so in this case we can use $\rho (T)$ instead of $\rho _{W}(Q_{\mathcal{P}%
},T) $ in all definitions we presented above.
\end{enumerate}
\end{remark}

\begin{remark}
If $T\in (Q_{\mathcal{P}}(X))_{0}$ has SVEP then for each $x\in X$ there
exists an unique maximal analytic extension of the application $\lambda
\rightarrow R(\lambda ,T)$, which will be denoted by $\tilde{x}$. Since $%
T\in (Q_{\mathcal{P}}(X))_{0}$ has SVEP the set $\rho _{T}(x)$ is
correctly defined and is unique. Moreover, $\rho _{T}(x)$ is open and $\sigma
_{T}(x)$ is closed.
\end{remark}

\begin{remark}
If $T\in (Q_{\mathcal{P}}(X))_{0}$ has SVEP and $x\in X$, then

\begin{enumerate}
\item $\rho _{T}(x)$ is an open set;

\item $\rho _{T}(x)$ is the domain of definition for $\tilde{x}$;

\item $\rho _{W}(Q_{\mathcal{P}},T)\subset \rho _{T}(x)$.
\end{enumerate}
\end{remark}

\begin{lemma}
{\label{lemma:qb1}} Let $(X,\mathcal{P})$ be a sequentially complete locally
convex space. If $T\in (Q_{\mathcal{P}}(X))_{0}$ then

\begin{enumerate}
\item the application $\lambda \rightarrow R(\lambda ,T)$ is holomorphic on $%
\rho _{W}(Q_{\mathcal{P}},T)$;

\item $\frac{d^{n}}{d\lambda ^{n}}R(\lambda ,T)=(-1)^{n}n!R(\lambda
,T)^{n+1} $, for every $n\in \mathbb{N}$;

\item $\lim\limits_{|\lambda |\rightarrow \infty }R(\lambda ,T)=0$ and $%
\lim\limits_{|\lambda |\rightarrow \infty }R(1,\lambda
^{-1}T)=\lim\limits_{|\lambda |\rightarrow \infty }\lambda R(1,T)=I$.
\end{enumerate}
\end{lemma}
\begin{proof}
1) If $\lambda _{0}\in \rho _{W}(Q_{\mathcal{P}},T)$ then there exists $V\in 
\mathcal{V}_{(\lambda _{0})}$ with the properties (1) and (2) from
definition of Walebroeck resolvent set. Since for every $\lambda \in
V\backslash \{\infty \}$ we have 
\begin{equation*}
R(\lambda ,T)-R(\lambda _{0},T)=(\lambda _{0}-\lambda )R(\lambda
,T)R(\lambda _{0},T)
\end{equation*}%
and the set $\{R(\lambda ,T)|~\lambda \in V\backslash \{\infty \}\}$ is
bounded in $Q_{\mathcal{P}}(X)$ results that the application $\lambda
\rightarrow R(\lambda ,T)$ is continuous in $\lambda _{0}$, so 
\begin{equation*}
\lim\limits_{\lambda \rightarrow \lambda_{0} }\frac{R(\lambda
,T)-R(\lambda_{0} ,T)}{\lambda -\lambda _{0}}=-R^{2}(\lambda _{0},T)
\end{equation*}

If $\lambda _{0}=\infty $ then, there exists some neighborhood $V\in 
\mathcal{V}_{(\infty )}$ such that the application $\lambda \rightarrow
R(\lambda ,T)$ is defined and bounded on $V\backslash \{\infty \}$.
Moreover, this application it is holomorphic and bounded on $V\backslash
\{\infty \}$, which implies that it is holomorphic at $\infty $.

Therefore, the application $\lambda \rightarrow R(\lambda ,T)$ is
holomorphic on $\rho _{W}(Q_{\mathcal{P}},T)$.

2) Results from the proof of (1).

3) For each $\lambda \in \rho _{W}(Q_{\mathcal{P}},T)$ we have 
\begin{equation*}
\lambda ^{-1}(I+TR(\lambda ,T))(\lambda I-T)=I,
\end{equation*}%
so 
\begin{equation}
{\label{equation:21}}R(\lambda ,T)=\lambda ^{-1}(I+TR(\lambda ,T)).
\end{equation}

If $V\in \mathcal{V}_{(\lambda _{0})}$ satisfies the conditions of the
definition of Walebroeck resolvent set, then the set 
\begin{equation*}
\{TR(\lambda ,T)|~\lambda \in V\backslash \{\infty \}\}
\end{equation*}%
is bounded, so from relation (\ref{equation:21}) results that $%
\lim\limits_{\vert\lambda\vert \rightarrow \infty }R(\lambda ,T)=0.$

From equality $R(\lambda ,T)=\lambda ^{-1}R(1,\lambda ^{-1}T),\lambda \neq 0$%
, and relation (\ref{equation:21}) results that 
\begin{equation*}
R(1,\lambda ^{-1}T)=I+TR(\lambda ,T),
\end{equation*}%
so 
\begin{equation*}
\lim\limits_{|\lambda |\rightarrow \infty }R(1,\lambda
^{-1}T)=\lim\limits_{|\lambda |\rightarrow \infty }(I+TR(\lambda ,T))=I
\end{equation*}
\end{proof}

\begin{lemma}
{\label{lemma:sv2}}If $T\in (Q_{\mathcal{P}}(X))_{0}$ has SVEP, then $\sigma
_{T}(x)=\varnothing $ if and only if $x=0_{X}$.
\end{lemma}
\begin{proof}
If $\sigma _{T}(x)=\varnothing $, then $\tilde{x}$ is an entire function.
Since $\left\vert \sigma _{W}(Q_{\mathcal{P}},T)\right\vert =r_{\mathcal{P}%
}(T)$, results that 
\begin{equation}
{\label{equation:ae1}}(\lambda I-T)\tilde{x}(\lambda )=x,~(\forall )|\lambda
|>r_{\mathcal{P}}(T),
\end{equation}%
so by lemma \ref{lemma:qb1} we have 
\begin{equation*}
\underset{|\lambda |\rightarrow \infty }{\lim }\tilde{x}(\lambda )=\underset{%
|\lambda |\rightarrow \infty }{\lim }R(\lambda ,T)x=0.
\end{equation*}

Therefore, from Liouville's theorem results that $\tilde{x}(\lambda )\equiv
0 $. Using the properties of functional calculus presented in \cite{st} and (\ref{equation:ae1}) we
have 
\begin{equation*}
x=\frac{1}{2\pi i}\int_{r_{\mathcal{P}}(T)+1}R(\lambda ,T)xd\lambda =\frac{1%
}{2\pi i}\int_{r_{\mathcal{P}}(T)+1}x(\lambda )d\lambda =0
\end{equation*}%
It is obvious that if $x=0_{X}$, then $\sigma _{T}(x)=\varnothing $.
\end{proof}

\section{ Quasi-nilpotent Equivalent Operators}

For a pair of operators $T,S\in (Q_{\mathcal{P}}(X))_{0}$, not necessarily
permutable, we consider the following notation

\begin{equation*}
(T-S)^{\left[ n\right] }=\sum%
\limits_{k=0}^{n}(-1)^{n-k}C_{n}^{k}T^{k}S^{n-k},
\end{equation*}%
where $C_{n}^{k}=\frac{n!}{\left( n-k\right) !k!}
$, for all $n\geq 1$ and $k=\overline{1,n}$.

\begin{remark}
\cite{co1} {\label{lemma:qn1}}If $T,S,P\in (Q_{\mathcal{P}}(X))_{0}$ then
for all $n\geq 1$ we have:

\begin{enumerate}
\item $(T-S)^{\left[ n+1\right] }=T(T-S)^{\left[ n\right] }-(T-S)^{\left[ n%
\right] }S.$

\item $\sum\limits_{k=0}^{n}(-1)^{n-k}C_{n}^{k}(T-S)^{\left[ k\right]
}(S-P)^{\left[ n-k\right] }=(T-P)^{\left[ n\right] }.$
\end{enumerate}
\end{remark}

\begin{definition}
We say that two operators $T,S\in (Q_{\mathcal{P}}(X))_{0}$ are
quasi-nilpotent equivalent operators if for every $p\in \mathcal{P}$ we have 
\begin{equation*}
\underset{n\rightarrow \infty }{\lim }\left( \hat{p}\left( (T-S)^{\left[ n%
\right] }\right) \right) ^{1/n}=0\text{ and }\underset{n\rightarrow \infty }{%
\lim }\left( \hat{p}\left( (T-S)^{\left[ n\right] }\right) \right) ^{1/n}=0.
\end{equation*}%
In this case we write $T\overset{q}{\backsim }$\ $S.$\ \ \ \ 
\end{definition}

\begin{remark}
If $T,S\in (Q_{\mathcal{P}}(X))_{0}$, then $(T-S)^{\left[ n\right] }\in Q_{%
\mathcal{P}}(X)$.
\end{remark}

\begin{lemma}
Let $(X,\mathcal{P})$ be a locally convex space and $T,S\in (Q_{\mathcal{P}%
}(X))_{0}$, such that $T\overset{q}{\backsim }S$. Then the series $%
\sum\limits_{n=0}^{\infty }(T-S)^{\left[ n\right] }$ and $%
\sum\limits_{n=0}^{\infty }(S-T)^{\left[ n\right] }$ converges in $Q_{%
\mathcal{P}}(X).$
\end{lemma}
\begin{proof}
If $T\overset{q}{\backsim }S$, then 
\begin{equation*}
\underset{n\rightarrow \infty }{\lim }\hat{p}\left( (T-S)^{\left[ n\right]
}\right) ^{1/n}=0,(\forall )p\in \mathcal{P},
\end{equation*}%
so by root test the series $\sum\limits_{n=0}^{\infty }\hat{p}((T-S)^{\left[
n\right] })$ converges. Moreover, for each $\varepsilon \in (0,1)$ and every 
$p\in \mathcal{P}$ there exists some index $n_{\varepsilon ,p}\in \mathbb{N}$
such that 
\begin{equation*}
\hat{p}\left( (T_{1}-T_{2})^{\left[ n\right] }\right) \leq \varepsilon
^{n},\left( \forall \right) n\geq n_{\varepsilon ,p}
\end{equation*}%
which implies that 
\begin{equation*}
\sum\limits_{k=n}^{m}\hat{p}\left( (T_{1}-T_{2})^{\left[ k\right] }\right)
<\sum\limits_{k=n}^{m}\varepsilon ^{k}<\frac{\varepsilon ^{n}}{1-\varepsilon 
},\left( \forall \right) \;m>n\geq n_{\varepsilon ,p},
\end{equation*}%
so $\left( \sum\limits_{k=0}^{n}(T_{1}-T_{2})^{\left[ k\right] }\right)
_{n\in \mathbb{N}}$ is a Cauchy sequence. Since the algebra $Q_{\mathcal{P}%
}(X)$ is sequentially complete results that the series $\sum\limits_{n=0}^{%
\infty }(T_{1}-T_{2})^{\left[ n\right] }$ converges in $Q_{\mathcal{P}}(X)$.

Analogously, we can prove that the series $\sum\limits_{n=0}^{\infty }(S-T)^{%
\left[ n\right] }$ converges in $Q_{\mathcal{P}}(X)$
\end{proof}

\begin{lemma}
{\label{lemma:qn2}}The relation $\overset{q}{\backsim }$ defined above is a
equivalence relation on $(Q_{\mathcal{P}}(X))_{0}$.
\end{lemma}
\begin{proof}
It is obvious that $\overset{q}{\backsim }$ is simetric and reflexive. Now
will prove that $\overset{q}{\backsim }$ is transitive. Let $%
T_{1},T_{2},T_{3}\in (Q_{\mathcal{P}}(X))_{0}$ such that $T_{1}\overset{q}{%
\backsim }T_{2}$ and $T_{2}\overset{q}{\backsim }T_{3}$. Then for every $%
\varepsilon >0$ and every $p\in \mathcal{P}$ there exists $n_{\varepsilon
,p}\in \mathbb{N}$ such that 
\begin{equation*}
\hat{p}\left( (T_{1}-T_{2})^{\left[ n\right] }\right) \leq \varepsilon ^{n}%
\text{ and }\hat{p}\left( (T_{2}-T_{3})^{\left[ n\right] }\right) \leq
\varepsilon ^{n},
\end{equation*}
for every $n\geq n_{\varepsilon ,p}$. If 
\begin{equation*}
M_{\varepsilon ,p}=\underset{k=\overline{1,n_{\varepsilon ,p}-1}}{\max }%
\left\{ \frac{\hat{p}\left( (T_{1}-T_{2})^{\left[ n\right] }\right) }{%
\varepsilon ^{k}},\frac{\hat{p}\left( (T_{2}-T_{3})^{\left[ n\right]
}\right) }{\varepsilon ^{k}},1\right\} ,(\forall )p\in \mathcal{P},
\end{equation*}%
then for every $n\in \mathbb{N}$ we have 
\begin{equation*}
\hat{p}\left( (T_{1}-T_{2})^{\left[ n\right] }\right) \leq M_{\varepsilon
,p}\varepsilon ^{n}\text{ and }\hat{p}\left( (T_{2}-T_{3})^{\left[ n\right]
}\right) \leq M_{\varepsilon ,p}\varepsilon ^{n},(\forall )p\in \mathcal{P}%
\text{ }.
\end{equation*}

The previous relation implies that 
\begin{equation*}
\hat{p}\left( (T_{1}-T_{3})^{\left[ n\right] }\right) =\hat{p}\left(
\sum\limits_{k=0}^{n}(-1)^{n-k}C_{n}^{k}(T_{1}-T_{2})^{\left[ k\right]
}(T_{2}-T_{3})^{\left[ n-k\right] }\right) \leq
\end{equation*}%
\begin{equation*}
\leq \sum\limits_{k=0}^{n}(-1)^{n-k}C_{n}^{k}\hat{p}\left( (T_{1}-T_{2})^{%
\left[ k\right] }\right) \hat{p}\left( (T_{2}-T_{3})^{\left[ n-k\right]
}\right)
\leq \sum\limits_{k=0}^{n}(-1)^{n-k}C_{n}^{k}M_{\varepsilon
,p}^{2}\varepsilon ^{k}\varepsilon ^{n-k}=\left( 2\varepsilon \right)^{n}M_{\varepsilon
,p}^{2}
\end{equation*}
for all $n\in \mathbb{N}$ and every $p\in \mathcal{P}$, so%
\begin{equation*}
\hat{p}\left( (T_{1}-T_{3})^{\left[ n\right] }\right)^{1/n} \leq 2\varepsilon \sqrt%
[n]{M_{\varepsilon,p }^{2}},\left( \forall \right) n\in \mathbb{N},(\forall
)p\in \mathcal{P}.
\end{equation*}

Therefore,%
\begin{equation*}
\underset{n\rightarrow \infty }{\lim }\hat{p}\left( (T_{1}-T_{3})^{\left[ n%
\right] }\right) ^{1/n}=0,(\forall )p\in \mathcal{P}
\end{equation*}

Analogously, we can prove that 
\begin{equation*}
\underset{n\rightarrow \infty }{\lim }\hat{p}\left( (T_{3}-T_{1})^{\left[ n%
\right] }\right) ^{1/n}=0,(\forall )p\in \mathcal{P},
\end{equation*}%
so $T_{1}\overset{q}{\backsim }$\ $T_{3}$.
\end{proof}

\begin{lemma}
If $(X,\mathcal{P})$ is a locally convex space then $T,S\in (Q_{\mathcal{P}%
}(X))_{0}$ are then quasi-nilpotent equivalent operators if and only if $%
\tilde{T}_{p},\tilde{S}_{p}\in \mathcal{L}(\tilde{X}^{p})$ are
quasi-nilpotent equivalent operators on the Banach space $\tilde{X}^{p}$,
for every $p\in \mathcal{P}$.
\end{lemma}
\begin{proof}
For every $p\in \mathcal{P}$ the subspace $N^{p}$ is invariant for $T_{p}$
and $T_{p}$, so 
\begin{equation*}
\left( \tilde{T}_{p}\right) ^{k}\left( \tilde{S}_{p}\right) ^{l}=\left(
T^{k}S^{l}\right) _{p},~(\forall )k,l\in \mathbb{N}.
\end{equation*}

Hence 
\begin{equation*}
(\tilde{T}_{p}-\tilde{S}_{p})^{\left[ n\right] }=\left( (T-S)^{\left[ n%
\right] }\right) _{p},~(\forall )m\in \mathbb{N}
\end{equation*}

If $T\overset{q}{\backsim }S$, then from definition results that 
\begin{equation}
\underset{n\rightarrow \infty }{\lim }\left( \hat{p}\left( (T-S)^{\left[ n%
\right] }\right) \right) ^{1/n}=0\text{ and }\underset{n\rightarrow \infty }{%
\lim }\left( \hat{p}\left( (T-S)^{\left[ n\right] }\right) \right) ^{1/n}=0.
\label{equation:qnil1}
\end{equation}%
so 
\begin{equation}
\underset{n\rightarrow \infty }{\lim }\left\Vert (\tilde{T}_{p}-\tilde{S}%
_{p})^{\left[ n\right] }\right\Vert _{p}^{1/n}=\underset{n\rightarrow \infty 
}{\lim }\hat{p}\left( (T-S)^{\left[ n\right] }\right) =0\text{ }
\label{equation:qnil2}
\end{equation}%
\begin{equation}
\underset{n\rightarrow \infty }{\lim }\left\Vert (\tilde{S}_{p}-\tilde{T}%
_{p})^{\left[ n\right] }\right\Vert _{p}^{1/n}=\underset{n\rightarrow \infty 
}{\lim }\hat{p}\left( (S-T)^{\left[ n\right] }\right) =0\text{ }
\label{equation:qnil3}
\end{equation}%
Therefore, $\tilde{T}_{p},\tilde{S}_{p}\in \mathcal{L}(\tilde{X}^{p}))_{0}$
are quasi-nilpotent equivalent operators, for every $p\in \mathcal{P}$.

Conversely, if $\tilde{T}_{p}\overset{q}{\backsim }\tilde{S}_{p}$, for every 
$p\in \mathcal{P}$, then the relation (\ref{equation:qnil2}) and (\ref%
{equation:qnil3}) holds, so condition (\ref{equation:qnil1}) is verified.
\end{proof}

\begin{lemma}
If $(X,\mathcal{P})$ is a locally convex space and $T,S\in (Q_{\mathcal{P}%
}(X))_{0}$ are then quasi-nilpotent equivalent operators, then $\sigma (Q_{%
\mathcal{P}},T)=\sigma (Q_{\mathcal{P}},S)$.
\end{lemma}

\begin{proof}
From previous lemma results that $\tilde{T}_{p},\tilde{S}_{p}\in \mathcal{L}(%
\tilde{X}^{p}))_{0}$ are quasi-nilpotent equivalent operators, for every $%
p\in \mathcal{P}$, hence by theorem 2.2 (\cite{co1}) we have $\sigma (%
\tilde{T}_{p})=\sigma (\tilde{S}_{p})$. Moreover, $\sigma (Q_{\mathcal{P}%
},T)=\underset{p}{\cup }\sigma (\tilde{T}_{p})$ and $\sigma (Q_{\mathcal{P}%
},S)=\underset{p}{\cup }\sigma (\tilde{S}_{p})$, so the corollary is proved.
\end{proof}

\begin{theorem}
Let $(X,\mathcal{P})$ be a locally convex space. If $T,S\in (Q_{\mathcal{P}%
}(X))_{0}$ are quasi-nilpotent equivalent operators, then $\sigma _{W}(Q_{%
\mathcal{P}},T)=\sigma _{W}(Q_{\mathcal{P}},S)$.
\end{theorem}

\begin{proof}
From lemma \ref{lemma:qb1} results that the functions $\lambda \rightarrow R(\lambda ,T)$ and 
$\lambda \rightarrow R(\lambda ,S)$ are holomorphic on the set $\rho _{W}(Q_{%
\mathcal{P}},T)$, respectively $\rho _{W}(Q_{\mathcal{P}},S)$.

Let $\lambda _{0}\in \sigma _{W}(Q_{\mathcal{P}},T)$ arbitrary fixed. Since $%
\sigma _{W}(Q_{\mathcal{P}},T)$ is an open set there exists $0<r_{1}<r_{2}$
such that $D_{i}(\lambda _{0})\subset \sigma _{W}(Q_{\mathcal{P}},T)$, $i=%
\overline{1,2}$, where%
\begin{equation*}
D_{i}(\lambda _{0})=\{\mu \in \mathbb{C}||\mu -\lambda _{0}|<r_{i}\},~i=%
\overline{1,2},
\end{equation*}%
and the set $\{~R(\lambda ,T)|~\lambda \in D_{1}(\lambda _{0})~\}$ is
bounded in $Q_{\mathcal{P}}(X)$. For each $p\in \mathcal{P}$ we consider
that 
\begin{equation*}
M_{p}=\sup \{~\hat{p}(R(\lambda ,T))|~\lambda \in D_{1}(\lambda _{0})~\}.
\end{equation*}

We denote by $R(\mu ,T)=\sum\limits_{k=0}^{n}R_{n}(\lambda )(\mu -\lambda
)^{n}$ the Taylor expansion of the resolvent around each point $\lambda $\
of $D_{1}(\lambda _{0})$. From complex analysis we have the formula%
\begin{equation*}
R_{n}(\lambda )=\frac{1}{n!}\frac{d^{n}}{d\lambda ^{n}}R(\mu ,T)=\frac{1}{%
2\pi i}\int\limits_{|\omega -\lambda |=r_{2}}\frac{R(\omega ,T)}{(\omega
-\lambda )^{n+1}}d\omega ,~(\forall )\lambda \in D_{1}(\lambda
_{0}),(\forall )n\geq 0,
\end{equation*}%
so,%
\begin{equation*}
\hat{p}(R_{n}(\lambda ))=\hat{p}\left( \frac{1}{2\pi i}\int\limits_{|\omega
-\lambda |=r_{2}}\frac{R(\omega ,T)}{(\omega -\lambda )^{n+1}}d\omega
\right) \leq
\end{equation*}%
\begin{equation*}
\leq \hat{p}(R_{n}(\lambda ))=r_{1}\sup \{~\hat{p}(R(\lambda ,T))|~\lambda
\in D_{1}(\lambda _{0})~\}\sup \{~\frac{1}{(\omega -\lambda )^{n+1}}%
|~\lambda \in D_{1}(\lambda _{0})~\}\leq
\end{equation*}%
\begin{equation*}
\leq r_{2}M_{p}(r_{2}-r_{1})^{-\left( n+1\right) },
\end{equation*}%
for all $\lambda \in D_{1}(\lambda _{0})$ and every $n\geq 0$. Since $T%
\overset{q}{\backsim }\ S$ results that for every $\varepsilon >0$ and
every $p\in \mathcal{P}$ there exists $n_{\varepsilon ,p}\in \mathbb{N}$
such that 
\begin{equation*}
\hat{p}\left( (T-S)^{\left[ n\right] }\right) \leq \varepsilon ^{n}\text{
and }\hat{p}\left( (S-T)^{\left[ n\right] }\right) \leq \varepsilon
^{n},(\forall )n\geq n_{\varepsilon ,p}.
\end{equation*}

Assume that $\varepsilon <r_{1}-r_{0}$. Then for every $p\in \mathcal{P}$
there exists $n_{\varepsilon ,p}\in \mathbb{N}$ such that
\begin{equation*}
\hat{p}\left( (S-T)^{\left[ n\right] }R_{n}(\lambda )\right) \leq
\varepsilon ^{n}r_{1}M_{p}(r_{1}-r_{0})^{-\left( n+1\right) }=
\end{equation*}%
\begin{equation*}
=r_{1}(r_{1}-r_{0})^{-1}M_{p}\left( \frac{\varepsilon }{r_{1}-r_{0}}\right) ^{n},
\end{equation*}
for every $n\geq n_{\varepsilon ,p}$ and every$~\lambda \in D_{1}(\lambda
_{0})$, so $\left( \sum\limits_{n=0}^{m}(-1)^{n}(S-T)^{\left[ n\right]
}R_{n}(\lambda )\right) _{m}$ is a Cauchy sequences. Since $Q_{\mathcal{P}%
}(X)$ is sequentially complete results that the series 
\begin{equation*}
R(\lambda
)=\sum\limits_{n=0}^{\infty }(-1)^{n}(S-T)^{\left[ n\right] }R_{n}(\lambda )
\end{equation*}
converges uniformly in $D_{0}$. Therefore, the function $\lambda \rightarrow 
$ $R(\lambda )$ is analytic in $\rho _{W}(Q_{\mathcal{P}},T)$.

Using lemma \ref{lemma:qb1} by induction it can be prove that if we
differentiate $n\geq 1$ times the equalities 
\begin{equation*}
(\lambda I-T)R(\lambda ,T)=R(\lambda ,T)(\lambda I-T)=I,
\end{equation*}%
then for every $n\geq 1$ we obtain 
\begin{equation*}
(\lambda I-T)\frac{d^{n}}{d\lambda ^{n}}R(\lambda ,T)=\frac{d^{n}}{d\lambda
^{n}}R(\lambda ,T)(\lambda I-T)=
\end{equation*}%
\begin{equation*}
=-n\frac{d^{n-1}}{d\lambda ^{n-1}}R(\lambda ,T),(\forall )\lambda \in
D_{0}(\lambda _{0})
\end{equation*}%
so%
\begin{equation*}
(\lambda I-T)R_{n}(\lambda )=(\lambda I-T)\frac{1}{n!}\frac{d^{n}}{d\lambda
^{n}}R(\mu ,T)=
\end{equation*}%
\begin{equation}
=-n\frac{1}{n!}\frac{d^{n-1}}{d\lambda ^{n-1}}R(\lambda ,T)=-R_{n-1}(\lambda
),  \label{equation:qn2}
\end{equation}%
for every $\lambda \in D_{1}(\lambda _{0})$ and every $n\geq 1$.

From lemma \ref{lemma:qn1} and relation (\ref{equation:qn2}) results the
following equalities 
\begin{equation*}
(\lambda I-S)R(\lambda )=\sum\limits_{n=0}^{\infty }(-1)^{n}(\lambda
I-S)(S-T)^{\left[ n\right] }R_{n}(\lambda )=
\end{equation*}
\begin{equation*}
=\sum\limits_{n=0}^{\infty }(\lambda I-S)\left( (\lambda I-S)-(\lambda
I-T)\right) ^{\left[ n\right] }SR_{n}(\lambda )=
\end{equation*}
\begin{equation*}
=\sum\limits_{n=0}^{\infty }\{ ((\lambda I-S)-(\lambda I-T))^{\left[ n+1%
\right] }R_{n}(\lambda )+((\lambda I-S)-(\lambda I-T))^{\left[ n\right]
}(\lambda I-T)R_{n}(\lambda )\}=
\end{equation*}
\begin{equation*}
=\sum\limits_{n=0}^{\infty }(-1)^{n+1}(S-T)^{\left[ n+1\right]
}R_{n}(\lambda )+(\lambda I-T)R_{0}(\lambda )+\sum\limits_{n=1}^{\infty
}(-1)^{n}(S-T)^{\left[ n\right] }(\lambda I-T)R_{n}(\lambda )=
\end{equation*}
\begin{equation*}
=\sum\limits_{n=0}^{\infty }(-1)^{n+1}(S-T)^{\left[ n+1\right]
}R_{n}(\lambda )+(\lambda I-T)R_{0}(\lambda )+\sum\limits_{n=0}^{\infty
}(-1)^{n+1}(S-T)^{\left[ n+1\right] }(\lambda I-T)R_{n+1}(\lambda )=
\end{equation*}
\begin{equation*}
=\sum\limits_{n=0}^{\infty }(-1)^{n+1}(S-T)^{\left[ n+1\right]
}R_{n}(\lambda )+(\lambda I-T)R(\lambda ,T)+\sum\limits_{n=0}^{\infty
}(-1)^{n+1}(S-T)^{\left[ n+1\right] }(-R_{n}(\lambda ))=I
\end{equation*}

Analogously we prove that $R(\lambda )(\lambda I-S)=I$, so $\rho _{W}(Q_{%
\mathcal{P}},T)\subset \rho _{W}(Q_{\mathcal{P}},S)$. The inclusion $\rho
_{W}(Q_{\mathcal{P}},S)\subset \rho _{W}(Q_{\mathcal{P}},T)$ can be proved
in the same way.
\end{proof}

\begin{lemma}
{\label{lemma:qb0}} Let $(X,\mathcal{P})$ be a locally convex space and $%
T\in (Q_{\mathcal{P}}(X))_{0}$ such that $r_{\mathcal{P}}(T)<1$. Then the
operator $I-T$ is invertible and $(I-T)^{-1}=\sum\limits_{n=0}^{\infty
}T^{n} $.
\end{lemma}
\begin{proof}
Assume that $r_{\mathcal{P}}(T)<t<1$. Hence results that 
\begin{equation*}
\underset{n\rightarrow \infty }{\limsup }\left( \hat{p}\left( T^{n}\right)
\right) ^{1/n}<t,\left( \forall \right) \;p\in \mathcal{P},
\end{equation*}
so for each $p\in \mathcal{P}$ there exists $n_{p}\in \mathbb{N}$ such that 
\begin{equation*}
(\hat{p}\left( T^{n}\right))^{1/n} \leq \underset{n\geq n_{p}}{\sup }\left( 
\hat{p}\left( T^{n}\right) \right) ^{1/n}<t,\left( \forall \right) \;n\geq
n_{p}.
\end{equation*}

This relation implies that the series $\sum\limits_{n=0}^{\infty }\hat{p}%
\left( T^{n}\right) $ converges, so 
\begin{equation*}
\underset{n\rightarrow \infty }{\lim }\hat{p}\left( T^{n}\right) =0,\left(
\forall \right) \;p\in \mathcal{P},
\end{equation*}%
therefore $\underset{n\rightarrow \infty }{\lim }T^{n}=0$. Since the algebra 
$Q_{\mathcal{P}}(X)$ is sequentially complete results that the series $%
\sum\limits_{n=0}^{\infty }T^{n}$ converges. Moreover, 
\begin{equation*}
(I-T)\sum\limits_{n=0}^{m}T^{n}=\sum\limits_{n=0}^{m}T^{n}(I-T)=I-T^{m+1},
\end{equation*}%
so 
\begin{equation*}
(I-T)\sum\limits_{n=0}^{\infty }T^{n}=\sum\limits_{n=0}^{\infty
}T^{n}(I-T)=I,
\end{equation*}%
which implies that $I-T$ is invertible and $(I-T)^{-1}=\sum\limits_{n=0}^{%
\infty }T^{n}$.
\end{proof}

\begin{theorem}
Let $(X,\mathcal{P})$ be a locally convex space. If $T,S\in (Q_{\mathcal{P}%
}(X))_{0}$ are quasi-nilpotent equivalent operators, then $T$ has SVEP if
and only if $S$ has SVEP.
\end{theorem}
\begin{proof}
Assume that $T$ has SVEP. Let $D_{f}\subset \mathbb{C}$ be an open set such
that $\rho _{W}(Q_{\mathcal{P}},S)\subset D_{f}$ and $f:D_{f}\rightarrow X$
be an analytic function on $D_{f}$ which satisfies the property%
\begin{equation*}
(\lambda I-S)f(\lambda )=0,~\lambda \in D_{f}.
\end{equation*}

Then, for every $n\geq 0$ we have%
\begin{equation*}
(T-S)^{\left[ n\right] }f(\lambda
)=\sum\limits_{k=0}^{n}(-1)^{n-k}C_{n}^{k}T^{k}S^{n-k}f(\lambda )=
\end{equation*}
\begin{equation}
=\sum\limits_{k=0}^{n}(-1)^{n-k}C_{n}^{k}T^{k}\lambda ^{n-k}f(\lambda
)=(T-\lambda I)^{n}f(\lambda )  \label{equation:qn3}
\end{equation}
Since $T\overset{q}{\backsim }\ S$ results that for every $\varepsilon >0$
and every $p\in \mathcal{P}$ there exists $n_{\varepsilon ,p}\in \mathbb{N}$
such that 
\begin{equation*}
\hat{p}\left( (T-S)^{\left[ n\right] }\right) \leq \varepsilon ^{n}\text{
and }\hat{p}\left( (S-T)^{\left[ n\right] }\right) \leq \varepsilon
^{n},(\forall )n\geq n_{\varepsilon ,p}.
\end{equation*}

Let $\mu \neq \lambda $. Then for every $\varepsilon \in (0,|\mu -\lambda |)$
and for every\ $p\in \mathcal{P}$ there exists $n_{\varepsilon ,p}\in 
\mathbb{N}$ 
\begin{equation*}
\hat{p}\left( \frac{(T-S)^{\left[ n\right] }}{\left( \mu -\lambda \right)
^{n+1}}\right) \leq \frac{\varepsilon ^{n}}{|\mu -\lambda |^{n+1}},(\forall
)n\geq n_{\varepsilon ,p},
\end{equation*}%
so the $\left( \sum\limits_{n=0}^{m}\frac{(T-S)^{\left[ n\right] }}{\left(
\mu -\lambda \right) ^{n+1}}\right) _{m}$ is a Cauchy sequences. Since $Q_{%
\mathcal{P}}(X)$ is sequentially complete results that the series $%
\sum\limits_{n=0}^{\infty }\frac{(T-S)^{\left[ n\right] }}{\left( \mu
-\lambda \right) ^{n+1}}$ is absolutely convergent in the topology of $Q_{%
\mathcal{P}}(X)$ for every $\mu \neq \lambda $. Moreover, if $r_{\mathcal{P}%
}(T-\lambda I)<|\mu -\lambda |$, then $r_{\mathcal{P}}(\frac{T-\lambda I}{%
\mu -\lambda })<1$ and from lemma \ref{lemma:qb0} results that 
\begin{equation}
\sum\limits_{n=0}^{\infty }\frac{(T-\lambda I)^{n}}{(\mu -\lambda )^{n+1}}%
=(\mu -\lambda )I-(T-\lambda I)=R(\mu ,T).  \label{equation:qn4}
\end{equation}%
From the relations (\ref{equation:qn3}) and (\ref{equation:qn4}) results
that 
\begin{equation*}
(\mu I-T)\left( \sum\limits_{n=0}^{\infty }\frac{(T-S)^{\left[ n\right] }}{%
(\mu -\lambda )^{n+1}}\right) f(\lambda )=(\mu I-T)\left(
\sum\limits_{n=0}^{\infty }\frac{(T-\lambda I)^{n}}{(\mu -\lambda )^{n+1}}%
\right) f(\lambda )=
\end{equation*}%
\begin{equation*}
=(\mu I-T)R(\mu ,T)f(\lambda )=f(\lambda ),
\end{equation*}%
for every $\mu $ with the property $r_{\mathcal{P}}(T-\lambda I)<|\mu
-\lambda |$. Therefore, 
\begin{equation*}
g_{\lambda }(\mu )=\sum\limits_{n=0}^{\infty }\frac{%
(T-S)^{\left[ n\right] }}{(\mu -\lambda )^{n+1}}f(\lambda ) 
\end{equation*}
is an analytic function on $\mathbb{C}\backslash \{\lambda \}$ which verifies the relation 
\begin{equation}
(\mu I-T)g_{\lambda }(\mu )=f(\lambda )  \label{equation:qn5}
\end{equation}%
on the open set $\{\mu \in \mathbb{C~}|r_{\mathcal{P}}(T-\lambda I)<|\mu
-\lambda |\}\subset \mathbb{C}\backslash \{\lambda \}$. Since T has SVEP
results that the function $g_{\lambda }(\mu )$ verifies the relation (\ref%
{equation:qn5}) for all $\mu \neq \lambda $. This implies that $\mathbb{C}%
\backslash \{\lambda \}\subset \rho _{T}(f(\lambda ))$, i.e. $\sigma
_{T}(f(\lambda ))\subset \{\lambda \}$.

Let $\lambda _{0}\in D_{f\text{ }}$ arbitrary fixed and $r>0$ such that $%
D_{0}=\{\lambda \in \mathbb{C~}|~|\lambda -\lambda _{0}|\leq r_{0}\}\subset
D_{f\text{ }}.$ Since $g_{\lambda }(\mu )$ is analytic on $\mathbb{C}%
\backslash \{\lambda \}$ from relation (\ref{equation:qn5}) results that 
\begin{equation}
(\mu I-T)\frac{1}{2\pi i}\int\limits_{|\xi -\lambda |=r_{0}}\frac{g_{\xi
}(\mu )}{\xi -\lambda _{0}}d\xi =\frac{1}{2\pi i}\int\limits_{|\xi -\lambda
|=r_{0}}\frac{f(\mu )}{\xi -\lambda _{0}}d\xi =f(\lambda _{0})
\label{equation:qn6}
\end{equation}%
for all $\mu \in D_{0}$, so $\mu \in \rho _{T}(f(\lambda _{0}))$, for every $%
\mu \in D_{0}$. Hence $\lambda _{0}\in \rho _{T}(f(\lambda _{0}))$ and since
we already proved above that $\sigma _{T}(f(\lambda _{0}))\subset \{\lambda _{0}\}$
results that $\sigma _{T}(f(\lambda ))=\dot{\varnothing}$. Lemma \ref%
{lemma:sv2} implies that $f(\lambda )\equiv 0$ on $D_{0}$ and since $\lambda
_{0}\in D_{f\text{ }}$ is arbitrary chosen, results that $f(\lambda )\equiv
0 $ on $D_{f}$. Therefore, $S$ has SVEP. Analogously we can prove that if $S$
has SVEP then $T$ has SVEP.
\end{proof}

\begin{theorem}
Let $(X,\mathcal{P})$ be a locally convex space. If $T,S\in (Q_{\mathcal{P}%
}(X))_{0}$ are quasi-nilpotent equivalent operators and $T$ has SVEP, then $%
\sigma _{T}(x)=\sigma _{S}(x)$, for every $x\in X$.
\end{theorem}
\begin{proof}
First we remark that from previous theorem results that $S$ has SVEP. Let $%
x(\lambda )$ be the analytic function on $\rho _{T}(x)$ which verify the
condition%
\begin{equation}
(\lambda I-T)x(\lambda )=x,~\lambda \in \rho _{T}(x).  \label{equation:qn7}
\end{equation}

Let $\lambda _{0}\in \rho _{T}(x)$ arbitrary fixed. Since $\rho _{T}(x)$ is
an open set there exists $0<r_{1}<r_{2}$ such that $D_{i}(\lambda
_{0})\subset \sigma _{W}(Q_{\mathcal{P}},T)$, $i=\overline{1,2}$, where%
\begin{equation*}
\bar{D}_{i}(\lambda _{0})=\{\mu \in \mathbb{C}|~|\mu -\lambda _{0}|\leq
r_{i}\},~i=\overline{1,2}.
\end{equation*}%
For every $\;p\in \mathcal{P}$ denote by $M_{p}^{1}$ the maximum of $%
x(\lambda )$\ on $\bar{D}_{2}(\lambda _{0})$. Hence, for $\lambda \in \bar{D}%
_{2}(\lambda _{0})$ we have 
\begin{equation}
p\left( \frac{x^{(n)}(\lambda )}{n!}\right) =p\left( \frac{1}{2\pi i}%
\int\limits_{|\xi -\lambda _{0}|=r_{2}}\frac{x(\xi )}{\left( \xi -\lambda
\right) ^{n+1}}d\xi \right) \leq \frac{M_{p}^{1}r_{2}}{(r_{2}-r_{1})^{n+1}},%
\text{ }(\forall )n\mathbb{\geq }0.  \label{equation:qn8}
\end{equation}%
In the proof of lemma \ref{lemma:qn2} we proved that for every $\varepsilon
>0$ and each $p\in \mathcal{P}$ there exists $M_{\varepsilon ,p}>0$ such
that 
\begin{equation}
\hat{p}\left( (T-S)^{\left[ n\right] }\right) \leq M_{\varepsilon
,p}\varepsilon ^{n},\text{ }(\forall )n\mathbb{\geq }0.  \label{equation:qn9}
\end{equation}%
Therefore, the relations (\ref{equation:qn8}) and (\ref{equation:qn9})
implies that 
\begin{equation*}
p\left( (T-S)^{\left[ n\right] }\frac{x^{(n)}(\lambda )}{n!}\right) \leq 
\hat{p}\left( (T-S)^{\left[ n\right] }\right) p\left( \frac{x^{(n)}(\lambda )%
}{n!}\right) <\frac{M_{\varepsilon ,p}M_{p}^{1}r_{2}}{r_{2}-r_{1}}\left( 
\frac{\varepsilon }{r_{2}-r_{1}}\right) ^{n},
\end{equation*}%
for all $n\mathbb{\geq }0$. Taking $\varepsilon =\frac{r_{2}-r_{1}}{2}$
results that for each $p\in \mathcal{P}$ there exists $M_{\varepsilon ,p}>0$
such that 
\begin{equation}
p\left( (T-S)^{\left[ n\right] }\frac{x^{(n)}(\lambda )}{n!}\right) \leq 
\frac{M_{p}}{2^{n}},\text{ }(\forall )n\mathbb{\geq }0,
\label{equation:qn10}
\end{equation}%
where $M_{p}=\frac{M_{\varepsilon ,p}M_{p}^{1}r_{2}}{r_{2}-r_{1}}$ does not
depend on $\lambda \in \bar{D}_{2}(\lambda _{0})$. The relation (\ref%
{equation:qn10}) shows that the series 
\begin{equation*}
\sum\limits_{n=0}^{\infty }p\left( \left( -1\right) ^{n}(T-S)^{\left[ n%
\right] }\frac{x^{(n)}(\lambda )}{n!}\right) 
\end{equation*}%
converges for every $\lambda \in \bar{D}_{2}(\lambda _{0})$ and every $p\in 
\mathcal{P}$, so since $X$ is sequentially complete results that the series $%
\sum\limits_{n=0}^{\infty }\left( -1\right) ^{n}(T-S)^{\left[ n\right] }%
\frac{x^{(n)}(\lambda )}{n!}$ converges absolutely and uniformly on $\bar{D}%
_{2}(\lambda _{0})$. But $\lambda _{0}\in \rho _{T}(x)$ is arbitrary fixed,
hence this series converges absolutely and uniformly on every compact $%
K\subset \rho _{T}(x)$, which implies that 
\begin{equation}
x_{1}(\lambda )=\sum\limits_{n=0}^{\infty }\left( -1\right) ^{n}(T-S)^{\left[
n\right] }\frac{x^{(n)}(\lambda )}{n!}  \label{equation:qn11}
\end{equation}%
is analytic on $\rho _{T}(x)$. Now we prove that 
\begin{equation*}
(\lambda I-S)x_{1}(\lambda )=x,~\lambda \in \rho _{T}(x).
\end{equation*}%
If we differentiate $n\geq 1$ times the equality (\ref{equation:qn7}), then
we have 
\begin{equation*}
(\lambda I-T)x^{(n)}(\lambda )=-nx^{(n-1)}(\lambda ),~\lambda \in \rho
_{T}(x).
\end{equation*}%
From previous relations and remark \ref{lemma:qn1} results 
\begin{equation*}
(\lambda I-S)x_{1}(\lambda )=\sum\limits_{n=0}^{\infty }(-1)^{n}(\lambda
I-S)(S-T)^{\left[ n\right] }\frac{x^{(n)}(\lambda )}{n!}=
\end{equation*}%
\begin{equation*}
=\sum\limits_{n=0}^{\infty }(\lambda I-S)\left( (\lambda I-S)-(\lambda
I-T)\right) ^{\left[ n\right] }\frac{x^{(n)}(\lambda )}{n!}=
\end{equation*}%
\begin{equation*}
=\sum\limits_{n=0}^{\infty }\left\{ \left( (\lambda I-S)-(\lambda
I-T)\right) ^{\left[ n+1\right] }+\left( (\lambda I-S)-(\lambda I-T)\right)
^{\left[ n\right] }(\lambda I-T)\right\} \frac{x^{(n)}(\lambda )}{n!}=
\end{equation*}%
\begin{equation*}
=\sum\limits_{n=0}^{\infty }(-1)^{n+1}\left( S-T\right) ^{\left[ n+1\right] }%
\frac{x^{(n)}(\lambda )}{n!}+\sum\limits_{n=0}^{\infty }(-1)^{n}\left(
S-T\right) ^{\left[ n\right] }(\lambda I-T)\frac{x^{(n)}(\lambda )}{n!}=
\end{equation*}%
\begin{equation*}
=\sum\limits_{n=0}^{\infty }(-1)^{n+1}\left( S-T\right) ^{\left[ n+1\right] }%
\frac{x^{(n)}(\lambda )}{n!}+(\lambda I-T)x(\lambda
)+\sum\limits_{n=1}^{\infty }(-1)^{n}\left( S-T\right) ^{\left[ n\right] }%
\frac{x^{(n-1)}(\lambda )}{\left( n-1\right) !}=
\end{equation*}%
\begin{equation*}
=(\lambda I-T)x(\lambda )=x
\end{equation*}
for all $\lambda \in \rho _{T}(x)$. This shows that $\rho _{T}(x)\subset
\rho _{S}(x)$, so $\sigma _{S}(x)\subset \sigma _{T}(x)$. Analogously it can
be proved that $\sigma _{T}(x)\subset \sigma _{S}(x)$.
\end{proof}

\begin{flushleft}
Sorin Mirel Stoian 

Department of Mathematics, University of Petro\c{s}ani

\textit{E-mail address:} mstoian@upet.ro
\end{flushleft}

\end{document}